\documentclass{amsart}

\newtheorem*{LocalDuality}{Local Duality}
\newtheorem{thm}{Theorem}
\newtheorem{cor}[thm]{Corollary}
\newtheorem{prop}[thm]{Proposition}
\newtheorem{defn}[thm]{Definition}
\newtheorem{lem}[thm]{Lemma}

\DeclareMathOperator{\depth}{depth}

\DeclareMathOperator{\Hom}{Hom}
\DeclareMathOperator{\mH}{H}
\DeclareMathOperator{\res}{res}
\newcommand{\m}{{\mathfrak m}}
\def\gfrac#1#2{\left[\begin{array}{c}#1\\#2\end{array}\right]}

\begin{document}


\title{Residues for Akizuki's one-dimensional local domain}
\author{I-Chiau Huang \and Jan-Li Lin}
\address{Institute of Mathematics, Academia Sinica, Nankang,
                 Taipei 11529, Taiwan, R.O.C. \and
                 Institute of Information Science, Academia Sinica, Nankang,
                 Taipei 11529, Taiwan, R.O.C.}
\email{ichuang@math.sinica.edu.tw \and ljl@iis.sinica.edu.tw}
\begin{abstract}
For a one-dimensional local domain $C_M$ constructed by Akizuki, we find
residue maps which give rise to a local duality.
The completion of $C_M$ is described using these residue maps.
\end{abstract}
\maketitle


Injective hulls of a given module are all isomorphic. For this reason, people
often speak of {\em the} injective hull to indicate its ``{\em uniqueness}''.
However, isomorphisms between these injective hulls are not canonical. In
fact, they are a part of the structure of the given module. For
instance, a local
duality for a power series ring \cite[(5.9)]{hu:pmzds} is interpreted
as an isomorphism between two injective hulls - one given by local cohomology
and another by continuous homomorphisms. This isomorphism is induced
by a residue map, which was not observed from the viewpoint of ``{\em
uniqueness}'' of injective hulls. In this article, our philosophy is
taken up again
by a Noetherian local ring $C_M$ constructed by Akizuki
\cite{aki:bpit}. Although
$C_M$ behaves beyond geometric expectation, we can still define certain
maps, which give rise to a local duality as an identification of
local cohomology
classes and continuous homomorphisms. These maps, also called residue maps,
determine all endomorphisms of an injective hull of the residue field of $C_M$.
So we are able to describe the completion of $C_M$.

We recall Akizuki's construction. Let $A$ be a discrete
valuation ring with the maximal ideal $\m=tA$, let $\hat{A}$ be its completion,
and let $K$ (resp. $\hat{K}$) be the quotient field of $A$ (resp. $\hat{A}$).
Assume that there is an element
\begin{eqnarray*}
z=a_0+a_1t^{n_1}+a_2t^{n_2}+\cdots\in\hat{A}&&(a_i\in A\setminus\m)
\end{eqnarray*}
transcendental over $A$ with the condition
\begin{eqnarray*}
n_r\geq2n_{r-1}+2&&(r\geq 1)
\end{eqnarray*}
on exponents, where $n_0=0$. Let
\begin{eqnarray*}
z_r=a_r+a_{r+1}t^{n_{r+1}-n_r}+\cdots&&(r\geq 0)
\end{eqnarray*}
and
$$
C=A[t(z_0-a_0),\{(z_i-a_i)^2\}_{i=0}^\infty].
$$
$C_M$ is defined to be the localization of $C$ at the maximal ideal $M$
generated by $t$ and $t(z_0-a_0)$. Akizuki \cite{aki:bpit} showed that $C_M$ is
a one-dimensional Noetherian local domain, whose normalization is not a finite
$C_M$-module.

The quotient field of $C_M$ is $K(z)$, which equals $(C_M)_t$ as $t$ is a
system of parameter of $C_M$. We can use the exact sequence
\begin{equation}
\label{eq:1}
0\to C_M\xrightarrow{\text{localization}}(C_M)_t\to\mH^1_{MC_M}(C_M)\to 0
\end{equation}
to describe elements of the first local cohomology module
$\mH^1_{MC_M}(C_M)$ of $C_M$ supported at the maximal ideal of $C_M$: For
$f\in C_M$ and $n>0$, the generalized fraction
$$
\gfrac{f}{t^n}_{C_M}
$$
is defined as the image of $f/t^n$ in $\mH^1_{MC_M}(C_M)$ under the map in
(\ref{eq:1}). The description of top local cohomology modules by generalized
fractions is essentially used in the concrete realizations of Grothendieck
duality \cite{hu:pmzds, hu:ecrc, hu:rtect}. It can be applied to combinatorial
analysis \cite{hu:arci,hu:rmca}. See also \cite{shar-zak:lcmgf} for
an alternate
treatment to generalized fractions.

For $f\in\hat{A}$ and $n>0$, we denote by
\begin{equation}
\label{eq:2}
\gfrac{f}{t^n}_A
\end{equation}
the image of $f/t^n$ in $\hat{K}/\hat{A}$ under the map in the exact sequence
$$
0\to\hat{A}\xrightarrow{\text{inclusion}}\hat{K}\to\hat{K}/\hat{A}\to 0.
$$
We remark that there is a canonical isomorphism between $\hat{K}/\hat{A}$ and
$\mH^1_{\m\hat{A}}(\hat{A})$, with which the representation (\ref{eq:2}) of
elements in $\hat{K}/\hat{A}$ and the representation of elements in
$\mH^1_{\m\hat{A}}(\hat{A})$ by generalized fractions agree.

Vanishing of elements in $\hat{K}/\hat{A}$ and $\mH^1_{MC_M}(C_M)$ can be
described in terms of ideal membership:
\begin{eqnarray*}
\gfrac{f}{t^n}_A=0
&&
\left(\text{ resp. $\gfrac{f}{t^n}_{C_M}=0$}\right)
\end{eqnarray*}
if and only if $f$ is contained in the ideal of $\hat{A}$ (resp.
$C_M$) generated
by $t^n$. For instance,
$$
\gfrac{t(z_0-a_0)}{t}_A=0
$$
but
$$
\gfrac{t(z_0-a_0)}{t}_{C_M}\neq 0.
$$
The canonical map
$K/A\to\hat{K}/\hat{A}$ is an isomorphism. We use the notation in (\ref{eq:2})
to represent elements in $K/A$ under this isomorphism.
\begin{lem}
\label{lem:1}
Every element of $\mH^1_{MC_M}(C_M)$ can be written as
$$
\gfrac{X+Yt(z_0-a_0)}{t^n}_{C_M}
$$
for some $X,Y\in A$, $n\ge 0$.
\end{lem}
\begin{proof}
It suffices to show that, for any specified $n> 0$, any element $f\in
C_M$ can be
written as
\begin{equation}
\label{eq:234}
f=X+Yt(z_0-a_0)+t^nZ
\end{equation}
with $X,Y\in A$ and $Z\in C_M$. Write $f$ as
$$
f=\frac{f_1}{1-f_2}
$$
for some $f_1\in C$ and $f_2\in M$. Then
$$
f-f_1(1+f_2+f_2^2+\cdots+f_2^n)=\frac{f_1f_2^{n+1}}{1-f_2}\in
M^{n+1}C_M.
$$
Since $M^{n+1}C_M\subset t^nC_M$,
$$
f-f_1(1+f_2+f_2^2+\cdots+f_2^n)=t^n Z_2
$$
for some $Z_2\in C_M$. By \cite[p. 140, Section 9.5, equation (5)]{rei:uca},
there exist $X,Y\in A$ and $Z_1\in C$ such that
$$
f_1(1+f_2+f_2^2+\cdots+f_2^n)=X+Yt(z_0-a_0)+t^nZ_1.
$$
The representation
$$
f=X+Yt(z_0-a_0)+t^n(Z_1+Z_2)
$$
is of the required form.
\end{proof}
\begin{lem}
\label{lem:2}
If $X,Y\in A$ and
$$
\gfrac{X+Yt(z_0-a_0)}{t^n}_{C_M}=0,
$$
then $X,Y\in t^n A$.
\end{lem}
\begin{proof}
$$
X+Yt(z_0-a_0)=t^nZ
$$
for some $Z\in C_M$. Write $X=t^\ell X_1$ and
$Y=t^m Y_1$ with invertible $X_1,Y_1\in A$. If $\ell\leq m$, then $n\leq\ell$,
otherwise
$X_1=-t^{m-\ell}Y_1t(z_0-a_0)+t^{n-\ell}Z\in t\hat{A}$. If $\ell>m$, then
$m\geq n$, otherwise
$t(z_0-a_0)=-t^{\ell-m}X_1Y_1^{-1}+t^{n-m}ZY_1^{-1}\in tC_M$. In either
case, $X,Y\in t^n A$.
\end{proof}
With these lemmas, we are able to define the following map for any
$\sigma,\rho\in\hat{A}$.
\begin{defn}
$$
\res_{\sigma,\rho}\colon\mH^1_{MC_M}(C_M)\to K/A
$$
is defined to be the $A$-linear map given by
$$
\gfrac{X+Yt(z_0-a_0)}{t^n}_{C_M}\mapsto
\gfrac{X\sigma+Y\rho}{t^n}_A.
$$
\end{defn}
Adopting the terminology of \cite{hu:pmzds}, we call $\res_{\sigma,\rho}$ a
residue map. Let
$$
\Hom_A^c(C_M,K/A)=
\{\varphi\in\Hom_A(C_M,K/A)\,|\,\varphi(M^nC_M)=0 \text{ for some $n$} \}
$$
be the $C_M$-module of continuous homomorphism. As a special case of J.
Lipman's result \cite[Proposition 3.4]{hu:pmzds}, $\Hom_A^c(C_M,K/A)$ is
an injective hull of the residue field of $C_M$. Note that
$t^{n+1}C_M\subset M^{n+1}C_M\subset t^nC_M$. Hence a $A$-linear map
$\varphi:C_M\to K/A$ is continuous if and only if $\varphi(t^nC_M)=0$ for some
$n$.  Using the representation (\ref{eq:234}) of elements of
$C_M$, we see that a continuous homomorphism is determined by an integer
$n$ with which $t^nC_M$ is in the kernel and by its values at $1$ and
$t(z_0-a_0)$.
\begin{defn}
$$
\Phi_{\sigma,\rho}\colon\mH^1_{MC_M}(C_M)\to\Hom_A^c(C_M,K/A)
$$
is defined to be the $C_M$-linear map given by
$$
\Phi_{\sigma,\rho}(\omega)(f)=\res_{\sigma,\rho}(f\omega),
$$
where $\omega\in\mH^1_{MC_M}(C_M)$ and $f\in C_M$.
\end{defn}
\begin{LocalDuality}
If $\rho$ is invertible, then $\Phi_{\sigma,\rho}$ is an isomorphism.
\end{LocalDuality}
\begin{proof}
The inverse map of $\Phi_{\sigma,\rho}$ can be written down explicitly. Let
$$
s_r:=a_1t^{n_1}+a_2t^{n_2}+\cdots+a_rt^{n_r}\in A.
$$
Then we have
$$
t(z_0-a_0)=t^{n_r+1}(z_r-a_r)+ts_r
$$
and
\begin{equation}
\label{eq:3}
t^2(z_0-a_0)^2+t^2s_r^2-2ts_rt(z_0-a_0)=t^{2n_r+2}(z_r-a_r)^2\in t^rC_M.
\end{equation}
Given $\varphi\in\Hom_A^c(C_M,K/A)$ with $\varphi(t^rC_M)=0$ and
$$\begin{cases}
\varphi(1)=\gfrac{\alpha}{t^r}_A\\
\varphi(t(z_0-a_0))=\gfrac{\beta}{t^r}_A,
\end{cases}$$
the system of equations
$$\begin{cases}
\gfrac{X\sigma+Y\rho}{t^r}_A=\gfrac{\alpha}{t^r}_A\\
\gfrac{X\rho-Yt^2s_r^2\sigma+2Yts_r\rho}{t^r}_A=\gfrac{\beta}{t^r}_A
\end{cases}$$
can be solved by choosing $X,Y\in A$ such that
$$\begin{cases}
X-\dfrac{\alpha ts_r(\sigma ts_r-2\rho)+\beta\rho}
              {(\rho-ts_r\sigma)^2}\in t^r\hat{A}
\\
Y-\dfrac{\alpha\rho-\beta\sigma}{(\rho-ts_r\sigma)^2}\in t^r\hat{A}.
\end{cases}$$
We define
$$
\Phi^{-1}(\varphi):=\gfrac{X+Yt(z_0-a_0)}{t^r}_{C_M},
$$
which is independent of the choices of $r$, $\alpha$, $\beta$, $X$ or $Y$. Then
$$
\Phi_{\sigma,\rho}(\Phi^{-1}(\varphi))(1)=\gfrac{\alpha}{t^r}_A
$$
and
\begin{eqnarray*}
& &
\Phi_{\sigma,\rho}(\Phi^{-1}(\varphi))(t(z_0-a_0))\\
&=&
\res_{\sigma,\rho}\gfrac{Xt(z_0-a_0)-Yt^2s_r^2+2Yts_rt(z_0-a_0)}{t^r}_{C_M}
\text{ (by (\ref{eq:3})) }\\
&=&
\gfrac{\beta}{t^r}_A.
\end{eqnarray*}
Hence $\Phi_{\sigma,\rho}(\Phi^{-1}(\varphi))=\varphi$. For any
$\omega\in\mH^1_{MC_M}(C_M)$, it is also straightforward to check that
$\Phi^{-1}(\Phi_{\sigma,\rho}(\omega))=\omega$. So $\Phi^{-1}$ is indeed the
inverse of $\Phi_{\sigma,\rho}$.
\end{proof}
\begin{cor}
$C_M$ is Gorenstein.
\end{cor}
\begin{proof}
$\mH^1_{MC_M}(C_M)$ is injective. Hence (\ref{eq:1}) is a finite injective
resolution of $C_M$.
\end{proof}
We remark that a Noetherian local ring $R$ whose maximal ideal is generated by
$1+\depth R$ elements is always Gorenstein \cite[p. 163, Exercise 1]{kap:cr}.
\begin{prop}
Any $C_M$-linear map
$$
\Phi\colon\mH^1_{MC_M}(C_M)\to\Hom_A^c(C_M,K/A)
$$
equals $\Phi_{\sigma,\rho}$ for some $\sigma,\rho\in\hat{A}$.
\end{prop}
\begin{proof}
For each $n$, there exist $\sigma_n,\rho_n\in A$ such that
$$\begin{cases}
\gfrac{\sigma_n}{t^n}_A=\Phi\left(\gfrac{1}{t^n}_{C_M}\right)(1)\\
\gfrac{\rho_n}{t^n}_A=\Phi\left(\gfrac{1}{t^n}_{C_M}\right)(t(z_0-a_0)).
\end{cases}$$
Since $\sigma_n-\sigma_{n+1}$ and $\rho_n-\rho_{n+1}$ are contained in
$t^nA$, the limits
$$\begin{cases}
\sigma=\lim_{n\to\infty}\sigma_n\\
\rho=\lim_{n\to\infty}\rho_n
\end{cases}$$
exist in $\hat{A}$. For any $X,Y\in A$,
$$
\Phi\left(\gfrac{X+Yt(z_0-a_0)}{t^n}_{C_M}\right)(1)=
\gfrac{X\sigma_n+Y\rho_n}{t^n}_A=
\gfrac{X\sigma+Y\rho}{t^n}_A.
$$
Hence
$$
(\Phi(\omega))(f)=\Phi(f\omega)(1)=
\Phi_{\sigma,\rho}(f\omega)(1)=(\Phi_{\sigma,\rho}(\omega))(f)
$$
for any $\omega\in\mH^1_{MC_M}(C_M)$ and $f\in C_M$. That is,
$\Phi=\Phi_{\sigma,\rho}$.
\end{proof}

Now we fix a $\sigma_0$ and an invertible $\rho_0$. All endomorphisms
of the $C_M$-module $\mH^1_{MC_M}(C_M)$ are of the form
$\Phi_{\sigma_0,\rho_0}^{-1}\circ\Phi_{\sigma,\rho}$. Since different pairs of
$\sigma$ and $\rho$ determine different $C_M$-linear maps
$\Phi_{\sigma,\rho}$, the completion $\widehat{C_M}$ of $C_M$ can be
described as the set
$\{\Phi_{\sigma,\rho}\}_{\sigma,\rho\in\hat{A}}$ with addition
$$
\Phi_{\sigma_1,\rho_1}+\Phi_{\sigma_2,\rho_2} =
\Phi_{\sigma_1+\sigma_2,\rho_1+\rho_2},
$$
unit $\Phi_{\sigma_0,\rho_0}$, and multiplication
$$
\Phi_{\sigma_1,\rho_1}*\Phi_{\sigma_2,\rho_2} =
\Phi_{\sigma_1,\rho_1}\circ
\Phi_{\sigma_0,\rho_0}^{-1}\circ\Phi_{\sigma_2,\rho_2}
$$
given by composition of endomorphisms. If $\sigma_0=0$ and $\rho_0=1$, then
$$
\Phi_{\sigma_1,\rho_1}*\Phi_{\sigma_2,\rho_2}=
\Phi_{\sigma_1\rho_2+\sigma_2\rho_1-2\sigma_1\sigma_2t(z_0-a_0),
\rho_1\rho_2-\sigma_1\sigma_2t^2(z_0-a_0)^2}.
$$
Identify $\Phi_{\sigma,\rho}$
with $\rho+\sigma X$ in $\hat{A}[X]/(X+t(z-a_0))^2$, compare their additions
and multiplications, we get the following description of $\widehat{C_M}$.
\begin{cor}
\label{cor:1}
$\widehat{C_M}\simeq\hat{A}[X]/(X+t(z-a_0))^2$
\end{cor}
For $\rho\in C_M$, the endomorphism $\Phi_{0,1}^{-1}\circ\Phi_{0,\rho}$ of
$\mH^1_{MC_M}(C_M)$ is multiplication by $\rho$. Therefore, with respect to
the isomorphism in Corollary~\ref{cor:1}, the embedding
$C_M\to\hat{A}[X]/(X+t(z-a_0))^2$ of completion is the composition
$$
C_M\xrightarrow{\text{inclusion}}\hat{A}\xrightarrow{\text{canonical}}
\hat{A}[X]/(X+t(z-a_0))^2.
$$


\end{document}